\documentclass[11pt]{article}
\usepackage{amsmath}
\usepackage{amssymb}
\usepackage{hyperref}
\usepackage{caption}
\usepackage{breakurl}
\usepackage{xurl}

\begin{document}

\title{Some Questions Around The Hilbert 16th Problem}
\author{Ali Taghavi \\ alitghv@yahoo.com\\}
\maketitle

Hilbert sixteenth problem asks for  existence of uniform upper
bound $H(n)$ for the number of limit cycles of a planar  polynomial vector
field of degree $n$. It is proved in \cite{Gaiko} that every polynomial vector field of degree $2$ has at most $4$ limit cycles.   It is well known that the number of limit cycles  of  a
polynomial vector field on the plane is finite, see \cite{Ilyashenko}. \\
Using polar change of coordinates we replace the concept of limit
cycle with $2\pi$-
periodic solution  $(r(t),\theta(t))$ for the corresponding vector field in $(r,\theta)$

More generally consider the equation\\
$$Z'=a_{n}(t)Z^{n}+a_{n-1}(t)Z^{n-1}+\ldots$$
where $Z'=dZ/dt$ and $a_{i}(t)$ are $c^{1}$ functions of real variable $t$.

For a fixed $t$, put $U(t)$ for the set of all $z$ in $\mathbb{C}$
with the property that  $\phi_{t}$ can be defined at $(Z,0)$, where $\phi$ is the
flow of corresponding autonomous system,adding $t'=1$ to the
equation. $U(t)$ is an open simply connected subset of
$\mathbb{C}$. Furthermore  $\phi_{t}$, as a map from $U(t)$ to
$\mathbb{C}$, is a holomorphic function. For if $a_{i}$ is not analytic
we approximate it with analytic functions by some Weierstrass type approximation theorem. Then we note that the uniform
limit of holomorphic functions is necessarily a  holomorphic function. This shows that for
$n>1$ the above system is not a complete Vector field. Because any
one to one entire function must be linear
 $\phi_{t}(Z)=a(t)Z+b(t)$. Differentiating in t implies that $n=1$
 Similarly if $\phi_{t}$ is a  mobious function then n can be at most two
This shows for $n>2$ the  argument of  Smale described in \cite{Sha} is
not applicable. Since every continuous function from sphere to sphere
is holomorphic if it is holomorphic on an open and dense subset of
spher  whose complement has zero Lesbegue measure. On the other hand every fixed point of a holomorphic function
has non negative Lefschetz number. In fact for$ n=3$ there are
examples of autonomous equation $z'=f(z)$ for which  there are two
isochronous centers with different period $T_{1}$ and $T_{2}$ so
$\phi_{T}$ can not be extended even continuously  to whole sphere.

The "\emph{complexification}" of the Hilbert 16th problem is an elegant and subtle
idea but in some cases is not effective.  In this note we suggest
some different  points for consideration of limit cycle problem.:

1)Let $[X,Y]=0$  and $\gamma$ be a limit cycle for X then $\gamma$
must be invariant under Y, namely $X$ and $Y$ share on limit
cycles. Since for every positive function $f$, the two vector fields $X$ and $fX$ have the same
trajectories, it is natural and interesting to compare   $C(X)$ with $C(fX)$.
By $C(X)$ , centralizer of $X$ ,one means all vector fields  $Z$
with $[X,Z]=0$. Note that locally around a non singular point of $X$, $C(X)$
and $C(fX)$ are isomorphic lie algebras. This local fact is no
longer true globally(for a non vanishing Vector field on arbitrary
surface),and it may be   false around a singularity of a vector field. For
example,put $f(x,y)=x^{2}+y^{2}+1$ and vector field $X$ as follow

 $$
\begin{cases}
 \dot{x}=y\\
 \dot{y}=-x,
\end{cases}
$$

$X$ is a non vanishing vector field on $\mathbb{R}^{2}\setminus\{0\}$, on the
other hand $X$ is a vector field with singularity at origin. In
both cases,two Lie algebras $C(X)$ and $C(fX)$ are not isomorph
since the operation of lie bracket is zero in $C(fX)$ ,but it is
not the case in $C(X)$.

It is also interesting that one look at the Hilbert sixteen
problem in a non analytic  but smooth  manner, for example
consider the following question:

Let $L$ be the Lienard polynomial vector field

$$
\begin{cases}
 \dot{x}=y-F(x)\\
 \dot{y}=-x,
\end{cases}
$$
without center and $S$ be a smooth vector field with $[L,S]=0$, is
it necessarily $S=kL$ for some constant $k$?\\

\textbf{ Remark} Non triviality of centralizer of non integrable
vector field X with components $(P,Q)$ is equivalent to complete
integrability of Hamiltonian $zP+wQ$ in $\mathbb{R}^{4}$
\newpage

\textbf{ Example 1} Consider vector fields $X$ and $Y$ as
follows:

$$(X)\;\;\;\;
\begin{cases}
 \dot{x}=y+x(x^{2}+y^{2}-1)\\
 \dot{y}=-x+y(x^{2}+y^{2}-1),
\end{cases}
$$

and

$$(Y)\;\;\;\;
\begin{cases}
 \dot{x}=2y+x(x^{2}+y^{2}-1)\\
 \dot{y}=-2x+y(x^{2}+y^{2}-1),
\end{cases}
$$

$X$ and $Y$ are independent out of circle $x^{2}+y^{2}=1$.This
circle is a hyperbolic limit cycle for X and Y while  we have
$[X,Y]=0$\\

2) We choose an arbitrary homogeneous vector field of degree
$n+1$.For example:

$$
\begin{cases}
 \dot{x}=y(x^{2}+y^{2})^{\frac{n}{2}}\\
 \dot{y}=-x(x^{2}+y^{2})^{\frac{n}{2}},
\end{cases}
$$

The cyclicity  of origin under perturbation of this vector field
among polynomial vector fields of  degree  at most $n$ is not less
than $H(n)$. Now put X for above Vector field  and Y for the
following linear center

$$
\begin{cases}
 \dot{x}=y\\
 \dot{y}=-x,
\end{cases}
$$
$[X,Y]=0$,Is it possible that for any perturbation $X_{\epsilon}$
of $X$,one have perturbation $Y_{\epsilon}$ with
$[X_{\epsilon},Y_{\epsilon}]=0$?

Note that linear center has finite $cyclicity$ by analytic
perturbation with finite parameter(The cyclicity of a singularity
of a vector field among a family of vector fields is the maximum
number of limit cycles which can be produced around the singularity  with
small perturbation of the vector field in the family,see \cite{Rouss}\\

 3)Let we have a vector field in $\mathbb{R}^{n}$ whose two first  components depend
 only on $x$ and $y$. Actually we have a planar  vector field
with these two components. Existence of invariant compact
$submanifold$ of codimension at most $2$ could lead to existence
of closed orbit in planar system. Since projection of such
submanifold on two first components  can not be a single point.\\

4)Any Vector field  on a surface defines a singular foliation of
dimension one and  according to definition  in \cite{X.Wang}  every limit cycle is considered as a  seperatrix.    It would be interesting to produce some
$C^{*}$ algebraic invariants depending only on degree of
polynomial Vector fields in the plane.\\

5)In \cite{FP} it is given a uniform upper bound, depending only on
$n$, for the length of a closed orbit of a polynomial vector field
of degree $n$. This computation is based on  the Riemannian metrics of the upper half sphere. On the other hand it
is well known that the number of closed geodesics of a surface
with negative curvature is uniformly bounded by the length of the
closed geodesics. Inspired by  these two concepts we ask: Can we equip the
phase portrait of a certain polynomial vector field  with an appropriate Riemannian metric whose curvature is
negative out of a finite number of analytic curves such that the
trajectories of our vector field would be unparameterized  geodesics. In particular is it possible to equip the punctured
plane to a Riemannian metric with negative curvature such that the
trajectories of the vander pol equation would be unparameterized geodesics. If the answer to the latter question is positive then
we could give another proof for the fact that the van der pol
equation can not have more than one limit cycle ,see \cite{hirsch-smale} for
information on the limit cycle of the vander pol equation. For a
subtle  relation between limit cycles and complex geometry ,see \cite{Romanovski}\\

6)A possible relation to operator theory: In the following arxiv note we
interpreted the number of limit cycles of the Lienard vector field
$L$ in terms of codimension of the range of functional operator
defined by $L$  with $L(g)=L.g$, the derivative of $g$ along the
trajectories of $L$. See \href{https://www.arxiv.org/abs/math.DS/0408037}{Counting limit Cycles via Index theory}
In fact if this operator happen to be bounded and has
a closed range  with respect to an appropriate norm, then we would have a    "Fredholm
index" interpretation for the number of limit cycles of our vector field. What Banach Functional space
is appropriate for the domain of the operator $L(g)=L.g$,such that
the operator would be a Fredholm operator whose  index is equal to the
number of limit cycles? Can we equip the space of smooth or
analytic maps on the plane with the structure of   a topological
vector space such that the corresponding operator would be a  Fredholm operator and the same
Fredholm index interpretation would be stile valid?(I thank A. Zeghib for his suggestion
for consideration of TVS as a possible resolution to this
problem). Is the generalization of the theory of "Fredholm
Operators on Banach space" for TVS a trivial problem?.

Finally we give two questions related to  subject 1:

\textbf{ Question 1 } For an analytic Vector field $X$ on the plane
or sphere  let  $C_{\omega}(X)$  be the space of all real analytic vector fields $Y$  with $[X,Y]=0$. Assume that $C_{\omega}(X)$ is an infinite dimensional Lie algebra. Does this imply that there is  a non constant analytic function $f$ globally defined on the plane or sphere with $X.f=0$?

\textbf{ Question 2 }.For every non vanishing vector field X on the
plane and positive smooth function f, is the centralizer $C(X)$ of $X$  isomorph to
$C(fX)$.If the  the answer is affirmative   we could actually assign  a
unique (up to isomorphism)  Lie algebra   to every smooth foliation of
the plane: We choose a generating vector field $X$ then we introduce $C(X)$ as the required Lie algebra.

\textbf{ Example 2}
Let an analytic vector field $X$ on sphere has a
center and a limit cycle simultaneously. Then its analytic centralizer $\mathcal{C}(X)$ consisting of all analytic vector field $Y$ with $[Y,X]=0$   has dimension one. Let  $V$ be the Van der pol
then $\mathcal{C}(V)$ has at most two dimension.(What is the exact dimension
of centralizer of Van der pol vector field?). In the example
1,$\mathcal{C}(X)$ is a 2 dimensional lie algebra while the centralizer of
the following vector field(as a vector field on the plane or
sphere) is a 4 dimensional space

$$
\begin{cases}
 \dot{x}=x\\
 \dot{y}=y,
\end{cases}
$$

\section*{Complex dilatation and Limit cycle theory}

In this current version of our note  we add this section  which is essentially the same as the materials in the following  Mathoveflow post: \\
 \href{https://mathoverflow.net/questions/382577/does-p-xp-yq-xq-y-0-implies-nonexistence-of-limit-cycle-for-p-partial-xq}{Limit cycle theory and  complex  dilitation}

Let $X=P\partial_x+Q\partial_y$ be a vector field on  $\mathbb{R}^2$. Assume that we have \begin{equation} \label{oui} P_xP_y+Q_xQ_y=0 \end{equation}
Does this imply that the vector field  $X$ is a divergence-free vector field  with respect to a Riemannian metric defined on $\mathbb{R}^2\setminus S$ where $S$ is the set  of singularities of $X$ ?
 Obviously this would imply that  any  $X$  with $P_xP_y+Q_xQ_y=0$ can not have any limit cycle. But is it really the case?
 Namely:

 \textbf{Question 3} Does the equality $P_xP_y+Q_xQ_y=0$ imply that $X=P\partial_x+Q\partial_y$  does not have any limit cycle?\\

 In the following we shall introduce some motivations for the above question:

 For  a  complex  function $f=U+iV:\mathbb{C} \to \mathbb{C}$ we recall the definition of the operators $f \mapsto \bar \partial f,\;\; f \mapsto  \partial f$  and the dilitation
  $\mu(f)$ as follows: $\bar \partial f=((U_x-V_y)+i(U_y+V_x),\;\;\;\partial f=(U_x+V_y)+i(V_x-U_y)\;\;\;\;\mu(f)=\frac{\bar \partial  f}{\partial f}$.

   To our vector field $X=P\partial_x+Q\partial_y$ we associate the complex function $X=P(z)+iQ(z)$.
 Then we consider the  complex  dilatation of $X$ with $\mu(z)=\bar\partial X(z)/\partial X(z)$.

 Let $\phi_t$ be the flow of the vector field $X$. So with a similar argument  for proof of   the standard variational equation in \cite[page 299]{hirsch-smale}  and  commutativity  of operatores  $\partial \phi_t$ and $\bar \partial \phi_t$  with $d/dt$ we  arrive at the following differential equations
 $$\begin{cases}(\partial \phi_t)'=\partial X(\phi_t(x))\partial \phi_t+\bar \partial X(\phi_t(x)) \partial \bar \phi_t\\ (\bar \partial \phi_t)'=\partial X(\phi_t(x)\bar \partial \phi_t+\bar \partial X(\phi_t(X))\overline{\partial g} \end{cases}$$
This imply that $\mu(\phi_t)'(x)=\frac{\bar\partial X(\phi_t(x))}{(\partial \phi_t)^2}exp\int_0^t div X(\phi_t (x))dt$\\

 \textbf{Proposition}\quad For a vector field $X=P\partial_x+Q\partial_y$ if the dilatation  $\mu(X(z))=\lambda $ is  a  constant map in $z$ with $|\lambda|=1$ then $X$ does not have any limit cycle.\\

\textbf{Proof}  Assume that $\mu(X(z))=\lambda$ for  a  fixed complex number $\lambda$ with $|\lambda|=1$ .  So  after a linear change of  coordinate $H(z)=(\beta )z $ where $\beta$ is a complex number with $\beta^2 \lambda=1$ then we have $H^* X(z)=\beta X(z/\beta)$. Since $H$ is a linear function then it is  equal to its linear part. Since both $H$  and $H^{-1}$ are holomorphic maps then we have  \begin{equation}\label{eqq}\begin{split} \mu(H^*(X))(z)=\mu (H^{-1}\circ X \circ H)(z)\\&=\mu(X)(z)\times (\frac{ H'(z)}{|H'(z)|})^2=\lambda.1/\beta^2=1\end{split}
   \end{equation}.

   For the  dilatation of functions  after right and left  composition with holomorphic maps  see \cite[page 9]{Ahlfors}

  Note that $X$ and $H^*(X)$ are orbitally equivalent vector fields. So if we prove that every vector field $Y$ with $\mu(Y(z))=1$ can not have any limit cycle then the proof of the proposition would be  completed. But it is an obvious fact: let $Y=P\partial_x+Q\partial_y$ be a vector field whose dilatation function is identically equal to $1$. this implies that $P_y=Q_y=0$. So $P,Q$ are functions in $x$. It is an standard fact in the theory of ordinary differential equation that a one dimensional autonomous vector field can not posses a periodic orbit. We apply this to $x'=P(x)$. So $Y$ has no periodic orbit. An alternative Proof for  non existence of periodic orbit for the planar vector field $Y(x,y)=P(x)\partial_x+Q(x)\partial_y$ is the following: We have $[Y,\partial/\partial y]=0$. If $\gamma$ is a periodic orbit of $Y$, then there is a point $A=(x_A,y_A)$ on $\gamma$ whose $x$-cordinate $x_A$ is maximum on $\gamma$. Then $Y(A)$ is a vertical vector,i.e. $Y(A)$ is parallel to the $y$-axis.So $P(A)=0$. Then $P(x_A,y)=0,\;\;\forall y\in \mathbb{R}$. This shows that the vertical line $x=x_A$ is invariant under the flow of $Y$. This contradicts the fact that the  periodic orbit $\gamma$ intersects this invariant line. In fact this situation violates the uniqueness of solutions of ordinary differential equation.

  So a natural question is the following:

  \textbf{Question 4} Is the above Proposition true without the assumption $|\lambda|=1 $?  According to the very process in ~\ref{eqq}  we may reformulate this question as follows; let  we have a vector field $X$  whose dilatation $\mu(X)=\lambda$ where $\lambda$ is a constant real number. Can such a vector field $X$ have any limit cycle? Obviously the columns of the jacobian of any vector field $X$ with the above property are orthogonal to each other. So this situation is a motivation for  question 3.

  \textbf{Some examples showing that questions 3 has possibly an
   affirmative answer;}

  Each of the following vector fields satisfy \eqref{oui} and non of them have any limit cycle:
  \enumerate
  \item
  A  vector field in the form $p(x)\partial_x+Q(y)\partial_y$ \\

  \item

  Every vector field in the form $P(y) \partial_x+Q(x) \partial_y$

  \item

  Every vector field in the form $Z'=f(z)$  where $f=P+iQ$ is a holomorphic map on $\mathbb{C}$

\end{document}